\documentclass[11pt,a4paper]{article}

\usepackage[margin=1in]{geometry}
\usepackage{amsmath,amssymb,amsthm}
\usepackage{mathtools}
\usepackage{mathrsfs}
\usepackage[T1]{fontenc}
\usepackage[utf8]{inputenc}
\usepackage[hidelinks]{hyperref}
\allowdisplaybreaks

\theoremstyle{plain}
\newtheorem{theorem}{Theorem}[section]
\newtheorem{lemma}[theorem]{Lemma}
\newtheorem{proposition}[theorem]{Proposition}
\newtheorem{corollary}[theorem]{Corollary}
\theoremstyle{definition}
\newtheorem{assumption}[theorem]{Assumption}
\theoremstyle{remark}
\newtheorem{remark}[theorem]{Remark}

\newcommand{\R}{\mathbb{R}}
\newcommand{\E}{\mathbb{E}}
\newcommand{\PP}{\mathbb{P}}
\newcommand{\Q}{\mathbb{Q}}
\newcommand{\F}{\mathcal{F}}
\newcommand{\Ecal}{\mathcal{E}}
\newcommand{\stexp}{\mathscr{E}}
\newcommand{\Law}{\mathcal{L}}
\newcommand{\dd}{\,d}

\DeclareMathOperator{\dist}{dist}

\title{Total variation bounds for Schr\"odinger--F\"ollmer sampling of
boundary-vanishing densities}

\author{Tatsuya Kaneko\textsuperscript{1}\qquad
Yumiharu Nakano\textsuperscript{1,2}\\[0.5em]
\small\textsuperscript{1}Department of Mathematical and Computing Science,
Institute of Science Tokyo\\
\small\textsuperscript{2}Future Intelligence Visionary Initiative,
Institute of Science Tokyo\\[0.4em]
\small\texttt{nakano@comp.isct.ac.jp}}

\date{November 2, 2026}

\begin{document}
\maketitle

\begin{abstract}
We derive total variation bounds for the Schr\"odinger--F\"ollmer sampler when
the target density vanishes at the boundary of its support. This boundary
behavior produces a terminal-time singularity in the drift near the boundary.
We first prove a base estimate under a single quantitative hypothesis,
finiteness of the relative Fisher information, without assuming boundedness or
Lipschitz continuity of the drift. We then combine this estimate with a mixture
regularization, thereby reducing the problem to determining the growth rate of
the regularized Fisher information. For densities with an exact polynomial
vanishing order on the boundary of a smooth bounded domain, we determine this
growth rate and its leading constant. Orders above one yield the square-root
rate in the number of steps, with a logarithmic loss at order one. For orders
between one half and one, the rate exponent decreases continuously from one
half toward one third. Finite mixtures of triangular kernel densities fall
into the linear case.
\end{abstract}

\noindent\textbf{Keywords.} Schr\"odinger--F\"ollmer sampler; Schr\"odinger
bridge; Euler--Maruyama scheme; total variation distance; relative Fisher
information.

\medskip
\noindent\textbf{MSC 2020.} 65C05; 65C30.

\section{Introduction and main results}\label{sec:intro}

The Schr\"odinger--F\"ollmer (SF) sampler reaches a prescribed target law at a
finite terminal time \cite{Follmer1986,DaiPra1991,HuangEtAl2025,DaiEtAl2023}.
Fix $T>0$, let
$G_T(x)=(2\pi T)^{-d/2}e^{-|x|^2/(2T)}$ and
$\gamma_T(dx)=G_T(x)\dd x$. Let $\mu(dx)=\rho(x)\dd x$ be the target
probability law and write $\phi=\rho/G_T$, so that $\mu=\phi\gamma_T$.
For $Z\sim N(0,I_d)$ define
\begin{equation}\label{eq:h-b}
h(t,x)=\E\bigl[\phi(x+\sqrt{T-t}\,Z)\bigr],
\qquad b(t,x)=\nabla\log h(t,x),\qquad 0\le t<T .
\end{equation}
The associated bridge
\cite{Schrodinger1931,Jamison1975,Follmer1986,DaiPra1991,Leonard2014} solves
\begin{equation}\label{eq:sf-sde}
dX_t=b(t,X_t)\dd t+dW_t,\qquad X_0=0,\qquad X_T\sim\mu,
\end{equation}
and its Euler--Maruyama approximation on $t_i=iT/n$ is
\begin{equation}\label{eq:euler-intro}
Y_{i+1}=Y_i+b(t_i,Y_i)\tfrac Tn+\sqrt{\tfrac Tn}\,Z_{i+1},\qquad Y_0=0 .
\end{equation}
Here and below, $(Z_i)$ are i.i.d.\ standard Gaussian random variables.

Existing quantitative estimates typically impose smoothness assumptions that
keep $b$ controlled up to time $T$
\cite{HuangEtAl2025,JiaoEtAl2026,DaiEtAl2023,WangZhang2026,LinWangZhang2026}.
This excludes densities with zeros, for which the drift is singular at the
terminal time. Endo and Nakano \cite{EndoNakano2024} addressed such targets by
the mixture regularization
\begin{equation}\label{eq:reg-intro}
\rho_\varepsilon=(1-\varepsilon)\rho+\varepsilon G_T,\qquad
\phi_\varepsilon=(1-\varepsilon)\phi+\varepsilon,\qquad
\mu_\varepsilon=\phi_\varepsilon\gamma_T ,
\end{equation}
but obtained only two-stage convergence, with no rate or admissible coupling
of $n$ and $\varepsilon$.

We first give the estimate that makes such a coupling possible. A Borel
function $f:\R^d\to[0,\infty)$ with $\int f\dd\gamma_T=1$ is called a density
ratio, and $\mu_f=f\gamma_T$. If $f$ has a weak gradient, set
\begin{equation}\label{eq:E-def}
\Ecal_T(f)=\int_{\R^d}\frac{|\nabla f|^2}{f}\dd\gamma_T\in[0,\infty],
\end{equation}
where the integrand is $0$ on $\{f=0,\nabla f=0\}$ and $+\infty$ on
$\{f=0,\nabla f\ne0\}$. For $f>0$, this is the relative Fisher information
$I(\mu_f\mid\gamma_T)$. Define $h,b$ and $Y_n$ for $f$ by replacing $\phi$
with $f$ in \eqref{eq:h-b}--\eqref{eq:euler-intro}. We use
$\delta(\nu_1,\nu_2)=\sup_{A\in\mathcal B(\R^d)}|\nu_1(A)-\nu_2(A)|$.

The next result controls the Euler--Maruyama discretization error by the
relative Fisher information.

\begin{theorem}[Base estimate]\label{thm:base}
Let $f:\R^d\to(0,\infty)$ be a density ratio such that
$f\in W^{1,1}_{\mathrm{loc}}(\R^d)$,
$\int e^{c|y|}\dd\mu_f(y)<\infty$ for some $c>0$, and
$\Ecal_T(f)<\infty$. Then, for every $n\ge1$,
\begin{equation}\label{eq:base}
\delta\bigl(\Law(Y_n),\mu_f\bigr)
\le\frac12\sqrt{\frac Tn\bigl(\Ecal_T(f)-|b(0,0)|^2\bigr)}
\le\frac12\sqrt{\frac Tn\,\Ecal_T(f)} .
\end{equation}
\end{theorem}

No boundedness, Lipschitz continuity, or time regularity of $b$ is required.
The proof derives an exact entropy identity and uses the martingale property
of the F\"ollmer drift; see \cite{EldanLehecShenfeld2020} for the latter.

For \eqref{eq:reg-intro}, let $b_\varepsilon$ be the drift associated with
$\phi_\varepsilon$, let $Y_n^\varepsilon$ be the corresponding output, and set
\begin{equation}\label{eq:Eeps}
\Ecal_\varepsilon:=\Ecal_T(\phi_\varepsilon)
=\int_{\R^d}\frac{(1-\varepsilon)^2|\nabla\phi|^2}
{(1-\varepsilon)\phi+\varepsilon}\dd\gamma_T .
\end{equation}

The next result separates the total variation error into discretization and
mixture regularization terms.

\begin{theorem}\label{thm:reg}
Let $\phi\in W^{1,1}_{\mathrm{loc}}(\R^d)$ and suppose that $\mu$ has a finite
exponential moment. If $\Ecal_\varepsilon<\infty$, then
\begin{equation}\label{eq:reg-bound}
\delta\bigl(\Law(Y^\varepsilon_n),\mu\bigr)
\le\frac12\sqrt{\frac Tn\,\Ecal_\varepsilon}+\varepsilon .
\end{equation}
\end{theorem}

Thus the problem is reduced to the boundary behavior of
$\Ecal_\varepsilon$. We impose the following precise condition.

\begin{assumption}\label{ass:boundary}
Let $\beta>0$. The probability density $\rho$ satisfies:
\begin{enumerate}
\renewcommand{\labelenumi}{(\alph{enumi})}
\item $\{\rho>0\}=D$ for a bounded open set $D\subset\R^d$ with $C^2$ boundary;
\item $\rho$ is continuous on $\R^d$ and locally Lipschitz in $D$;
\item writing $\nu$ for the inner unit normal, there is a continuous
$c:\partial D\to(0,\infty)$ such that, uniformly in $y\in\partial D$ as
$r\downarrow0$,
\[
\rho\bigl(y+r\nu(y)\bigr)=c(y)r^\beta(1+o(1)),\qquad
\nabla\rho\bigl(y+r\nu(y)\bigr)
=\beta c(y)r^{\beta-1}\bigl(\nu(y)+o(1)\bigr).
\]
\end{enumerate}
\end{assumption}

Let $\sigma$ denote surface measure on $\partial D$.

The next result shows how the boundary vanishing order determines the
regularized Fisher information and the resulting convergence rate.

\begin{theorem}\label{thm:main}
Under Assumption~\ref{ass:boundary}, the following statements hold.
\begin{enumerate}
\renewcommand{\labelenumi}{(\roman{enumi})}
\item If $\beta>1$, then $\Ecal_0:=\Ecal_T(\phi)<\infty$ and
$\Ecal_\varepsilon\le\Ecal_0$. With $\varepsilon_n=n^{-1/2}$,
\[
\delta\bigl(\Law(Y_n^{\varepsilon_n}),\mu\bigr)
\le\frac12\sqrt{\frac Tn\,\Ecal_0}+\frac1{\sqrt n}=O(n^{-1/2}).
\]
\item If $\beta=1$, then
\[
\lim_{\varepsilon\downarrow0}
\frac{\Ecal_\varepsilon}{\log(1/\varepsilon)}
=\int_{\partial D}c\dd\sigma .
\]
With $\varepsilon_n=n^{-1/2}$,
$\delta(\Law(Y_n^{\varepsilon_n}),\mu)
=O(\sqrt{\log n/n})$.
\item If $\tfrac12<\beta<1$, then
\[
\lim_{\varepsilon\downarrow0}
\varepsilon^{(1-\beta)/\beta}\Ecal_\varepsilon
=\frac{\pi\beta}{|\sin(\pi/\beta)|}
\int_{\partial D}G_T^{1-1/\beta}c^{1/\beta}\dd\sigma
=:\Lambda_\beta\in(0,\infty).
\]
In particular, with $\varepsilon_n\asymp n^{-\beta/(1+\beta)}$,
$\delta(\Law(Y_n^{\varepsilon_n}),\mu)
=O(n^{-\beta/(1+\beta)})$.
\item If $\beta\le\tfrac12$, then $\Ecal_\varepsilon=\infty$ for every
$\varepsilon\in(0,1)$.
\end{enumerate}
\end{theorem}

For $\beta=1$, the Gaussian factor cancels and the leading coefficient is the
total inward boundary slope. The range $\frac12<\beta<1$ is not covered by the
growth-weighted Lipschitz condition of \cite{EndoNakano2024}; our base estimate
applies because it does not require a bounded drift. Section~\ref{sec:reg}
gives a generic $O(n^{-1/3})$ bound, while Section~\ref{sec:boundary} proves
the sharp energy asymptotics above and treats triangular kernel mixtures
explicitly.

\section{Proof of the base estimate}\label{sec:base}

Let $\Omega=C([0,T],\R^d)$ with coordinate process $X_t(\omega)=\omega(t)$, let
$\PP$ be the Wiener measure starting at $0$, and let $(\F_t)_{0\le t\le T}$ be the
usual augmentation of the raw filtration. For a density ratio $f$ define
$\PP^f$ on $\F_T$ by $d\PP^f/d\PP=f(X_T)$; this is a probability measure and
$\PP^f\circ X_T^{-1}=\mu_f$. Expectations under $\PP$ and $\PP^f$ are written
$\E$ and $\E^f$. Write
$p_s(x,y)=(2\pi s)^{-d/2}e^{-|y-x|^2/(2s)}$ for the Brownian transition density,
and define $h$ and $b$ by \eqref{eq:h-b}, i.e.\
$h(t,x)=\int f(y)p_{T-t}(x,y)\dd y$ and $b=\nabla\log h$.

For $n\ge1$ put $\Delta_n=T/n$, $t_i=i\Delta_n$ and
$\tau_n(t)=t_i$ for $t\in[t_i,t_{i+1})$, $i=0,\dots,n-1$, with
$\tau_n(T):=t_{n-1}$ (the value at the single point $T$ is immaterial). On an
auxiliary atomless complete probability space $(\Omega_E,\mathcal G,\Q)$
carrying i.i.d.\ $N(0,I_d)$ variables $(Z_i)_{i\ge1}$, the $n$-step SF sampler
is defined by \eqref{eq:euler-intro}. For any random variable $Y$ on this
space, write $\Law(Y)=\Q\circ Y^{-1}$.

We compare path measures using
$H(Q\mid R)=\int\log(dQ/dR)\dd Q$ if $Q\ll R$, and $+\infty$ otherwise.

\subsection{The F\"ollmer drift process}

\begin{lemma}\label{lem:follmer}
Under the hypotheses of Theorem~{\rm \ref{thm:base}}, the following hold.
\begin{enumerate}
\renewcommand{\labelenumi}{\textup{(\alph{enumi})}}
\item $h(t,x)\in(0,\infty)$ for all $t\in(0,T)$, $x\in\R^d$, and also for $t=0$
and $|x|<cT$; on this set $h$ is smooth in $x$, satisfies
$\partial_th+\frac12\Delta h=0$, and
\begin{equation}\label{eq:grad-h}
\nabla h(t,x)=\int_{\R^d}\frac{y-x}{T-t}\,f(y)p_{T-t}(x,y)\dd y
=\int_{\R^d}\nabla f(y)\,p_{T-t}(x,y)\dd y .
\end{equation}
In particular $h(0,0)=1$ and $b(0,0)=\int\nabla f\dd\gamma_T$ is finite.
\item Setting $v_t=b(t,X_t)$ for $0\le t<T$, one has
$v_t=\E^f[\nabla\log f(X_T)\mid\F_t]$; consequently $v$ extends to a
square-integrable $\PP^f$-martingale on $[0,T]$ with
\begin{equation}\label{eq:vT}
v_T=\nabla\log f(X_T),\qquad \E^f|v_T|^2=\Ecal_T(f),\qquad
v_0=b(0,0),
\end{equation}
and $\E^f\int_0^T|v_t|^2\dd t\le T\,\Ecal_T(f)<\infty$.
\item Under $\PP^f$ the process $W^f_t=X_t-\int_0^tb(s,X_s)\dd s$,
$t\in[0,T]$, is a Brownian motion; thus $X$ solves \eqref{eq:sf-sde} under
$\PP^f$ and $X_T\sim\mu_f$.
\end{enumerate}
\end{lemma}

\begin{proof}
\textup{(a).} Fix $t\in(0,T)$ and set $s=T-t\in(0,T)$. Then
\begin{equation}\label{eq:kernel-domination}
\frac{p_s(x,y)}{G_T(y)}
=\Bigl(\frac Ts\Bigr)^{d/2}
\exp\Bigl(-\frac12\Bigl[\Bigl(\frac1s-\frac1T\Bigr)|y|^2
-\frac{2x\cdot y}{s}+\frac{|x|^2}{s}\Bigr]\Bigr),
\end{equation}
and since $s<T$, the coefficient of $|y|^2$ in the exponent is strictly
negative. Hence, for each $(t,x)$, there are $C=C(t,x)$ and $a=a(t)>0$ such that
$p_{T-t}(x,y)\le C e^{-a|y|^2}G_T(y)$, locally uniformly in $(t,x)$. Therefore
$h(t,x)\le C\int f\dd\gamma_T=C<\infty$, and the same domination justifies
differentiating the kernel under the integral sign, which gives the first
equality in \eqref{eq:grad-h} and, by iteration, smoothness in $x$; the backward
heat equation follows because $p_{T-t}(x,y)$ solves it in $(t,x)$. Positivity of
$h$ is clear from $f>0$.

For $t=0$ we have $p_T(0,y)=G_T(y)$, so $h(0,0)=\int f\dd\gamma_T=1$, and more
generally
$p_T(x,y)=G_T(y)e^{x\cdot y/T-|x|^2/(2T)}$, and therefore
$h(0,x)\le e^{-|x|^2/(2T)}\int e^{|x||y|/T}\dd\mu_f(y)<\infty$ for $|x|<cT$ by
the exponential-moment hypothesis of Theorem~\ref{thm:base}; the same bound
dominates the difference quotients
and yields \eqref{eq:grad-h} at $t=0$ for $|x|<cT$.

For the second equality in \eqref{eq:grad-h}, note first that by
Cauchy--Schwarz and \eqref{eq:kernel-domination},
\begin{equation}\label{eq:gradf-int}
\begin{aligned}
\int|\nabla f(y)|p_{T-t}(x,y)\dd y
&\le C\int|\nabla f|\dd\gamma_T \\
&\le C\Bigl(\int\tfrac{|\nabla f|^2}{f}\dd\gamma_T\Bigr)^{1/2}
\Bigl(\int f\dd\gamma_T\Bigr)^{1/2}
=C\sqrt{\Ecal_T(f)}<\infty .
\end{aligned}
\end{equation}
Let $\chi_R$ be smooth cut-offs with $\mathbf 1_{B_R}\le\chi_R\le\mathbf
1_{B_{2R}}$ and $|\nabla\chi_R|\le 2/R$. Since $f\in W^{1,1}_{\mathrm{loc}}$,
\[
\int\nabla f(y)\chi_R(y)p_{T-t}(x,y)\dd y
=-\int f(y)\bigl[\nabla\chi_R(y)p_{T-t}(x,y)
+\chi_R(y)\nabla_yp_{T-t}(x,y)\bigr]\dd y .
\]
The first term is bounded by $(2/R)\int f\,p_{T-t}(x,\cdot)\to0$, and
$\nabla_yp_{T-t}(x,y)=-\frac{y-x}{T-t}p_{T-t}(x,y)$; letting $R\to\infty$ and
using \eqref{eq:gradf-int} together with the finiteness of
$\int|y-x|f(y)p_{T-t}(x,y)\dd y$ (again \eqref{eq:kernel-domination}) gives
\eqref{eq:grad-h}. Finally $b(0,0)=\nabla h(0,0)/h(0,0)=\int\nabla f\dd\gamma_T$.

\textup{(b).} Let $0\le t<T$ (with $t=0$ admissible since $X_0=0$). By the Markov
property, $\E[f(X_T)\mid\F_t]=h(t,X_t)$ and
$\E[\nabla f(X_T)\mid\F_t]=\int\nabla f(y)p_{T-t}(X_t,y)\dd y$, both finite a.s.
by (a) and \eqref{eq:gradf-int}. Bayes' formula and \eqref{eq:grad-h} give
\[
\E^f\bigl[\nabla\log f(X_T)\mid\F_t\bigr]
=\frac{\E[f(X_T)\nabla\log f(X_T)\mid\F_t]}{\E[f(X_T)\mid\F_t]}
=\frac{\E[\nabla f(X_T)\mid\F_t]}{\E[f(X_T)\mid\F_t]}
=\frac{\nabla h(t,X_t)}{h(t,X_t)}=v_t .
\]
Since $X_T\sim\mu_f$ under $\PP^f$,
$\E^f|\nabla\log f(X_T)|^2=\int|\nabla\log f|^2\dd\mu_f=\Ecal_T(f)<\infty$, so
$(v_t)_{t<T}$ is bounded in $L^2$ and therefore extends to a
square-integrable martingale with terminal value $\nabla\log f(X_T)$.
Moreover, $v_0=b(0,0)$ because $\F_0$ is trivial. By Jensen's inequality,
$\E^f|v_t|^2\le\E^f|v_T|^2$ for all $t$, and hence
$\E^f\int_0^T|v_t|^2\dd t\le T\Ecal_T(f)$.

\textup{(c).} Let $M_t=\E[f(X_T)\mid\F_t]$, a uniformly integrable $\PP$-martingale
with $M_0=1$, $M_T=f(X_T)$ and $M_t=h(t,X_t)>0$ for $t<T$. By It\^o's formula
and the backward heat equation, $dM_t=\nabla h(t,X_t)\cdot dX_t=M_tb(t,X_t)\cdot
dX_t$ on $(0,T)$, and by continuity of $M$ at $0$ this representation holds on
$[0,T)$. Hence $M$ is the stochastic exponential of $\int_0^\cdot b(s,X_s)\cdot
dX_s$ on $[0,T)$, and $M_T=f(X_T)>0$, so $\PP^f$ and $\PP$ are equivalent on
$\F_T$. By (b), $\E^f\int_0^T|v_t|^2\dd t<\infty$, in particular
$\int_0^T|b(s,X_s)|\dd s<\infty$ $\PP^f$-a.s.\ and, by equivalence, $\PP$-a.s.
The Girsanov--Maruyama theorem for the density process $M$ therefore applies
on all of $[0,T]$: $W^f=X-\int_0^\cdot b(s,X_s)\dd s$ is a continuous $\PP^f$-local
martingale with bracket $\langle W^{f,i},W^{f,j}\rangle_t=\delta_{ij}t$, hence a
Brownian motion by L\'evy's characterization. Finally, for Borel $A$,
$\PP^f(X_T\in A)=\E[f(X_T)\mathbf 1_A(X_T)]=\int_Af\dd\gamma_T=\mu_f(A)$.
\end{proof}

\subsection{Proof of Theorem~\ref{thm:base}}

\begin{proof}
Step 1.
Put $\beta_t^{(n)}=b(\tau_n(t),X_{\tau_n(t)})$, an adapted process that is
constant on each $[t_i,t_{i+1})$; all its values are finite by
Lemma~\ref{lem:follmer}(a), because $\tau_n(t)\le t_{n-1}<T$ and $X_0=0$. Hence
$\int_0^T|\beta^{(n)}_t|^2\dd t=\Delta_n\sum_{i<n}|b(t_i,X_{t_i})|^2<\infty$
a.s., and
\[
Z_t=\exp\Bigl(\int_0^t\beta^{(n)}_s\cdot dX_s
-\tfrac12\int_0^t|\beta^{(n)}_s|^2\dd s\Bigr)
\]
is a positive continuous $\PP$-local martingale. It is a true martingale: on
$[t_i,t_{i+1}]$ the integrand is the $\F_{t_i}$-measurable constant vector
$c_i=b(t_i,X_{t_i})$, so
$\E\bigl[e^{c_i\cdot(X_{t_{i+1}}-X_{t_i})-\frac12|c_i|^2\Delta_n}\mid
\F_{t_i}\bigr]=1$, and iterating over $i=n-1,\dots,0$ gives $\E Z_T=1$. (No
Novikov-type condition, hence no boundedness of $b$, is needed here.)

Define $\PP_n$ by $d\PP_n/d\PP=Z_T$; then $\PP_n\sim\PP$ and, by the
Girsanov--Maruyama theorem,
$W^{(n)}=X-\int_0^\cdot\beta^{(n)}_s\dd s$ is a $\PP_n$-Brownian motion.
Consequently, for $i=0,\dots,n-1$,
\[
X_{t_{i+1}}=X_{t_i}+b(t_i,X_{t_i})\Delta_n+\sqrt{\Delta_n}\,\xi^{(n)}_{i+1},
\qquad
\xi^{(n)}_{i+1}=\frac{W^{(n)}_{t_{i+1}}-W^{(n)}_{t_i}}{\sqrt{\Delta_n}},
\]
with $(\xi^{(n)}_i)_{i=1}^n$ i.i.d.\ $N(0,I_d)$ under $\PP_n$. Comparing with
\eqref{eq:euler-intro},
\begin{equation}\label{eq:grid-law}
\PP_n\circ(X_{t_0},\dots,X_{t_n})^{-1}=\Q\circ(Y_0,\dots,Y_n)^{-1},
\qquad\text{in particular}\quad \PP_n\circ X_T^{-1}=\Law(Y_n).
\end{equation}

Step 2.
By Lemma~\ref{lem:follmer}(c) and Step 1, $\PP^f$ and $\PP_n$ are both
equivalent to $\PP$. Recall from Lemma~\ref{lem:follmer}(b) that
$v_t=b(t,X_t)$ is the F\"ollmer drift process of the exact bridge under
$\PP^f$. With $\stexp$ denoting the Dol\'eans--Dade stochastic exponential,
\[
\frac{d\PP^f}{d\PP}=\stexp\left(\int v\cdot dX\right)_T,
\qquad
\frac{d\PP_n}{d\PP}
=\stexp\left(\int\beta^{(n)}\cdot dX\right)_T.
\]
The following calculation is first performed after stopping the stochastic
integrals when their quadratic variations reach an integer level. Substituting
$dX=dW^{(n)}+\beta^{(n)}dt$ gives
\[
\log\frac{d\PP^f}{d\PP_n}
=\int_0^T(v_t-\beta^{(n)}_t)\cdot dW^{(n)}_t
-\frac12\int_0^T|v_t-\beta^{(n)}_t|^2\dd t .
\]
Under $\PP^f$ we have $dW^{(n)}=dW^f+(v-\beta^{(n)})\dd t$, so
\[
\log\frac{d\PP^f}{d\PP_n}
=\int_0^T(v_t-\beta^{(n)}_t)\cdot dW^f_t
+\frac12\int_0^T|v_t-\beta^{(n)}_t|^2\dd t .
\]
By Lemma~\ref{lem:follmer}(b) and Jensen,
\begin{equation}\label{eq:L2-drift}
\E^f\int_0^T|v_t-\beta^{(n)}_t|^2\dd t
\le2\,\E^f\!\int_0^T\!|v_t|^2\dd t
+2\Delta_n\sum_{i=0}^{n-1}\E^f|v_{t_i}|^2
\le4T\,\Ecal_T(f)<\infty,
\end{equation}
so the stopped stochastic integrals converge in $L^2(\PP^f)$ to the stochastic
integral against $W^f$, which is a true $\PP^f$-martingale with zero mean.
Letting the stopping levels tend to infinity in the log-density identity then
gives
\begin{equation}\label{eq:entropy}
H(\PP^f\mid\PP_n)
=\frac12\,\E^f\int_0^T|v_t-v_{\tau_n(t)}|^2\dd t .
\end{equation}

Step 3.
Let $t\in[t_i,t_{i+1}]$. Since $v$ is a square-integrable martingale, its
increments over $[t_i,t]$ and $[t,t_{i+1}]$ are orthogonal in $L^2(\PP^f)$, so
$\E^f|v_{t_{i+1}}-v_{t_i}|^2=\E^f|v_t-v_{t_i}|^2+\E^f|v_{t_{i+1}}-v_t|^2$ and
hence $\E^f|v_t-v_{t_i}|^2\le\E^f|v_{t_{i+1}}-v_{t_i}|^2$. Summing and using
orthogonality once more,
\[
\E^f\int_0^T|v_t-v_{\tau_n(t)}|^2\dd t
\le\Delta_n\sum_{i=0}^{n-1}\E^f|v_{t_{i+1}}-v_{t_i}|^2
=\Delta_n\,\E^f|v_T-v_0|^2
=\Delta_n\bigl(\Ecal_T(f)-|b(0,0)|^2\bigr),
\]
the last equality because $v_0=\E^f[v_T]=b(0,0)$ is deterministic. With
\eqref{eq:entropy}, $H(\PP^f\mid\PP_n)\le\frac12\Delta_n(\Ecal_T(f)-|b(0,0)|^2)$.

Step 4.
The map $\pi_T(\omega)=\omega(T)$ pushes $\PP^f$ to $\mu_f$ and, by
\eqref{eq:grid-law}, $\PP_n$ to $\Law(Y_n)$. Relative entropy does not increase
under measurable maps, so
$H(\mu_f\mid\Law(Y_n))\le H(\PP^f\mid\PP_n)$, and Pinsker's inequality gives
\[
\delta(\Law(Y_n),\mu_f)\le\sqrt{\tfrac12H(\mu_f\mid\Law(Y_n))}
\le\frac12\sqrt{\Delta_n\bigl(\Ecal_T(f)-|b(0,0)|^2\bigr)} . \qedhere
\]
\end{proof}

\begin{remark}\label{rem:complexity}
The quantity controlling the discretization cost is more precisely
\[
\mathcal I_T(f):=\Ecal_T(f)-|b(0,0)|^2
=\E^f|v_T-v_0|^2.
\]
Thus, for a prescribed accuracy $\eta\in(0,1)$ in total variation, it
suffices to choose
\[
n\ge
\max\left\{1,\left\lceil\frac{T\mathcal I_T(f)}{4\eta^2}\right\rceil\right\}.
\]
Here $\lceil\cdot\rceil$ denotes the ceiling function.
Since $\mathcal I_T(f)\le\Ecal_T(f)$, the simpler condition
\[
n\ge\max\left\{1,\left\lceil
\frac{T\Ecal_T(f)}{4\eta^2}\right\rceil\right\}
\]
is also sufficient.
In particular, for fixed $T$ and $f$, the required number of time steps is
$O(\eta^{-2})$, with no boundedness or Lipschitz constant of the drift entering
the estimate. If $\mathcal I_T(f)=0$, the bound in
Theorem~\ref{thm:base} vanishes for every $n\ge1$; in this case $v$ is
constant and the Euler scheme is exact at the terminal time.
\end{remark}

\section{Regularization bounds}\label{sec:reg}

\begin{proof}[Proof of Theorem~{\rm \ref{thm:reg}}]
$\phi_\varepsilon\ge\varepsilon>0$ and $\int\phi_\varepsilon\dd\gamma_T=1$;
$\nabla\phi_\varepsilon=(1-\varepsilon)\nabla\phi$ so
$\phi_\varepsilon\in W^{1,1}_{\mathrm{loc}}$; and
$\mu_\varepsilon=(1-\varepsilon)\mu+\varepsilon\gamma_T$ has a finite
exponential moment since both $\mu$ and $\gamma_T$ do. Thus the hypotheses of
Theorem~\ref{thm:base} hold for $f=\phi_\varepsilon$, and
Theorem~\ref{thm:base} gives
$\delta(\Law(Y^\varepsilon_n),\mu_\varepsilon)\le\frac12\sqrt{T\Ecal_\varepsilon/n}$.
Since
$\delta(\mu_\varepsilon,\mu)=\varepsilon\,\delta(\gamma_T,\mu)\le\varepsilon$,
the triangle inequality completes the proof.
\end{proof}

The problem is thus reduced to estimating $\Ecal_\varepsilon$. The following
generic bound uses the condition introduced in \cite{EndoNakano2024}.

\begin{assumption}\label{ass:growth}
There is $C_1>0$ such that
$|\phi(x)-\phi(y)|\le C_1(1+\phi(x)+\phi(y))|x-y|$ for all $x,y\in\R^d$.
\end{assumption}

The next result derives a generic $O(n^{-1/3})$ convergence rate from this
growth condition.

\begin{proposition}\label{prop:generic}
Under Assumption~{\rm \ref{ass:growth}}, the hypotheses of
Theorem~{\rm \ref{thm:reg}} hold and, for every $\varepsilon\in(0,1]$,
\begin{equation}\label{eq:Eeps-generic}
\Ecal_\varepsilon\le C_1^2\bigl(\varepsilon^{-1}+2-3\varepsilon\bigr)
\le\tfrac43C_1^2\varepsilon^{-1}. 
\end{equation}
Moreover, with
$\varepsilon_n=\min\{1,(C_1\sqrt{T/3}/2)^{2/3}n^{-1/3}\}$,
\[
\delta\bigl(\Law(Y^{\varepsilon_n}_n),\mu\bigr)=O(n^{-1/3}).
\]
\end{proposition}

\begin{proof}
The case $\varepsilon=1$ ($\phi_\varepsilon\equiv1$) is trivial, so let
$\varepsilon\in(0,1)$. Fix $x_0$ and $r>0$ with
$C_1r\le\frac12$. For $|x-x_0|\le r$,
Assumption~\ref{ass:growth} gives
$(1-C_1r)\phi(x)\le(1+C_1r)\phi(x_0)+C_1r$, which implies
$\phi(x)\le3\phi(x_0)+1$. Iterating along a chain of
$\lceil2C_1|x|\rceil$ balls of radius $r=1/(2C_1)$ yields
$\phi(x)\le(1+\phi(0))3^{\,2C_1|x|+1}$. Since $\gamma_T$ is Gaussian, $\mu$
therefore has exponential moments of every order.
Local boundedness of $\phi$ combined with Assumption~\ref{ass:growth} shows that
$\phi$ is locally Lipschitz, hence in $W^{1,1}_{\mathrm{loc}}$, and Rademacher's
theorem gives $|\nabla\phi(x)|\le C_1(1+2\phi(x))$ a.e.

Therefore, with $r=\phi(x)$ and $D=(1-\varepsilon)r+\varepsilon$, the algebraic
identity
\[
\frac{(1-\varepsilon)^2(1+2r)^2}{D}
=4(1-\varepsilon)r+4(1-2\varepsilon)+\frac{(1-3\varepsilon)^2}{D}
\qquad(r\ge0)
\]
(verified by multiplying out; the constant term matches because
$4\varepsilon(1-2\varepsilon)+(1-3\varepsilon)^2=(1-\varepsilon)^2$) together
with $\int\phi\dd\gamma_T=1$ and $D\ge\varepsilon$ gives
\[
\Ecal_\varepsilon\le C_1^2\int\frac{(1-\varepsilon)^2(1+2\phi)^2}{D}\dd\gamma_T
\le C_1^2\Bigl(4(1-\varepsilon)+4(1-2\varepsilon)
+\frac{(1-3\varepsilon)^2}{\varepsilon}\Bigr)
=C_1^2\bigl(\varepsilon^{-1}+2-3\varepsilon\bigr).
\]
Finally,
$\varepsilon^{-1}+2-3\varepsilon
=\varepsilon^{-1}(1+2\varepsilon-3\varepsilon^2)\le\frac4{3\varepsilon}$,
because the maximum of $1+2\varepsilon-3\varepsilon^2$ is $\frac43$, attained at
$\varepsilon=\frac13$. Substituting into \eqref{eq:reg-bound} gives
$\delta\le A n^{-1/2}\varepsilon^{-1/2}+\varepsilon$ with $A=C_1\sqrt{T/3}$,
minimized at $\varepsilon=(A/2)^{2/3}n^{-1/3}$.
\end{proof}

\begin{remark}\label{rem:q}
In the proof of Proposition~\ref{prop:generic}, we estimate
$\Ecal_\varepsilon$ directly from \eqref{eq:Eeps}, keeping the factor
$((1-\varepsilon)\phi+\varepsilon)^{-1}$ inside the integral. This gives the
bound $\Ecal_\varepsilon=O(\varepsilon^{-1})$. A cruder argument would first
use the following uniform estimate. Since
$\nabla\phi_\varepsilon=(1-\varepsilon)\nabla\phi$ and
$|\nabla\phi|\le C_1(1+2\phi)$ a.e.,
\[
|\nabla\log\phi_\varepsilon|
\le C_1\left(
\frac{1-\varepsilon}{(1-\varepsilon)\phi+\varepsilon}
+\frac{2(1-\varepsilon)\phi}{(1-\varepsilon)\phi+\varepsilon}
\right)
\le C_1\left(\frac{1-\varepsilon}{\varepsilon}+2\right).
\]
Integrating the square of this estimate against $\mu_\varepsilon$ would give only
$\Ecal_\varepsilon=O(\varepsilon^{-2})$. Balancing the two terms in
\eqref{eq:reg-bound} would then yield the rate $n^{-1/4}$ instead of the
$n^{-1/3}$ rate in Proposition~\ref{prop:generic}.
More generally, if $\Ecal_\varepsilon\le C\varepsilon^{-q}$ for some $q>0$, then
\eqref{eq:reg-bound} with $\varepsilon_n\asymp n^{-1/(q+2)}$ gives
$\delta(\Law(Y^{\varepsilon_n}_n),\mu)=O(n^{-1/(q+2)})$.
Proposition~\ref{prop:generic} is the case $q=1$.
\end{remark}

\section{Boundary asymptotics}\label{sec:boundary}

We now prove Theorem~\ref{thm:main}. We first explain the role of the boundary
condition in Assumption~\ref{ass:boundary}. Compact support without controlled
vanishing at the boundary does not ensure that the regularized Fisher
information is finite. For example, if $\rho=\mathbf 1_{[0,1]}$, then
$\phi_\varepsilon\notin W^{1,1}_{\mathrm{loc}}$ and
$\Ecal_\varepsilon=\infty$ for every $\varepsilon$. The mixture
\eqref{eq:reg-intro} makes $\phi_\varepsilon$ positive but does not remove the
jump at the boundary. This explains why the boundary condition is needed.

For $\tfrac12<\beta<1$, define
\begin{equation}\label{eq:Ibeta}
I(\beta):=\int_0^\infty\frac{u^{2\beta-2}}{1+u^\beta}\dd u
=\frac{\pi}{\beta\,|\sin(\pi/\beta)|}\in(0,\infty).
\end{equation}

The exponent in (iii) increases continuously from $\frac13$ to $\frac12$ as
$\beta$ ranges over $(\frac12,1)$. At $\beta=1$, it matches the rate in (ii)
up to the logarithmic factor; at the excluded endpoint $\beta=\frac12$, it
matches the generic exponent $\frac13$ in Proposition~\ref{prop:generic}.
Thus the exponent $q=1$ in
Remark~\ref{rem:q} is attained exactly in the limit $\beta\downarrow\frac12$,
where $\Ecal_\varepsilon$ ceases to be finite. All four regimes are governed by
the boundary alone, and the critical exponent $\frac12$ does not depend on the
dimension.

\subsection{A one-sided boundary lemma}

The following one-dimensional lemma describes the contribution to
$\Ecal_\varepsilon$ along one inward normal direction at the boundary.

\begin{lemma}\label{lem:boundary}
Let $\beta>0$, $\kappa>0$, $g>0$, and $\delta>0$. Let
$\psi:(0,\delta)\to(0,\infty)$ be continuous, and let
$a:(0,\delta)\to[0,\infty)$ and $w:(0,\delta)\to(0,\infty)$ be measurable.
Suppose that, for every $r_0\in(0,\delta)$,
\[
\inf_{r_0\le r<\delta}\psi(r)>0,
\qquad
\int_{r_0}^{\delta}a(r)^2w(r)\dd r<\infty,
\]
and that, as $r\downarrow0$,
\[
\psi(r)=\kappa r^\beta(1+o(1)),\quad
a(r)=\beta\kappa r^{\beta-1}(1+o(1)),\quad
w(r)=g(1+o(1)).
\]
Define
$J_\varepsilon=\int_0^\delta\dfrac{(1-\varepsilon)^2a(r)^2w(r)}
{(1-\varepsilon)\psi(r)+\varepsilon}\dd r$. Then:
\begin{enumerate}
\renewcommand{\labelenumi}{\textup{(\roman{enumi})}}
\item if $\beta>1$, $\displaystyle\sup_{0<\varepsilon\le1}J_\varepsilon
\le\int_0^\delta\frac{a^2w}{\psi}\dd r<\infty$;
\item if $\beta=1$, $\displaystyle\lim_{\varepsilon\downarrow0}
\frac{J_\varepsilon}{\log(1/\varepsilon)}=\kappa g$;
\item if $\tfrac12<\beta<1$, $\displaystyle\lim_{\varepsilon\downarrow0}
\varepsilon^{(1-\beta)/\beta}J_\varepsilon=\beta^2\,g\,\kappa^{1/\beta}I(\beta)$;
\item if $\beta\le\tfrac12$, $J_\varepsilon=\infty$ for every
$\varepsilon\in(0,1)$.
\end{enumerate}
For the uniform version, let $\mathcal A$ be an index set and attach a subscript
$\alpha$ to $\psi,a,w,\kappa,g$, and $J$. Suppose that $\beta$ and $\delta$ are
fixed, the three asymptotic relations above hold uniformly in
$\alpha\in\mathcal A$, the set
$\{(\kappa_\alpha,g_\alpha):\alpha\in\mathcal A\}$ is contained in a compact
subset of $(0,\infty)^2$, and, for every $r_0\in(0,\delta)$,
\[
\inf_{\alpha\in\mathcal A}\inf_{r_0\le r<\delta}\psi_\alpha(r)>0,
\qquad
\sup_{\alpha\in\mathcal A}\int_{r_0}^{\delta}
a_\alpha(r)^2w_\alpha(r)\dd r<\infty.
\]
Then, in case (i),
$\sup_{\alpha\in\mathcal A}\sup_{0<\varepsilon\le1}J_{\varepsilon,\alpha}<\infty$;
in cases (ii) and (iii), respectively,
\[
\sup_{\alpha\in\mathcal A}\left|
\frac{J_{\varepsilon,\alpha}}{\log(1/\varepsilon)}
-\kappa_\alpha g_\alpha\right|\longrightarrow0
\qquad(\varepsilon\downarrow0),
\]
\[
\sup_{\alpha\in\mathcal A}\left|
\varepsilon^{(1-\beta)/\beta}J_{\varepsilon,\alpha}
-\beta^2g_\alpha\kappa_\alpha^{1/\beta}I(\beta)\right|\longrightarrow0
\qquad(\varepsilon\downarrow0).
\]
\end{lemma}

\begin{proof}
(i) Since $(1-\varepsilon)\psi+\varepsilon\ge(1-\varepsilon)\psi$, the integrand
is bounded by $a^2w/\psi$, which near $0$ behaves like
$\beta^2\kappa g\,r^{\beta-2}$ and is therefore integrable for $\beta>1$; on
$[\delta/2,\delta)$ its integral is finite by hypothesis and the positive lower
bound on $\psi$.

(iv) Near $0$ the integrand is at least
$(1-\varepsilon)^2a^2w/(\varepsilon+\psi)$, which has exact order
$r^{2\beta-2}$ because $\psi(r)\to0$; moreover,
$\int_0^{\delta_0}r^{2\beta-2}\dd r=\infty$ for every $\delta_0>0$ when
$2\beta-2\le-1$.

(ii)--(iii) Let $\eta\in(0,1)$ and choose $\delta_\eta\in(0,\delta)$ such that on
$(0,\delta_\eta)$,
\[
(1-\eta)\kappa r^\beta\le\psi(r)\le(1+\eta)\kappa r^\beta,\qquad
(1-\eta)g\le w(r)\le(1+\eta)g,
\]
\[
(1-\eta)\beta\kappa r^{\beta-1}\le a(r)\le(1+\eta)\beta\kappa r^{\beta-1}.
\]
Set
$m_\eta=\inf_{\delta_\eta\le r<\delta}\psi(r)>0$. Since
\[
(1-\varepsilon)\psi(r)+\varepsilon
\ge\min\{m_\eta,1\},
\qquad \delta_\eta\le r<\delta,
\]
and $(1-\varepsilon)^2\le1$, the contribution of
$[\delta_\eta,\delta)$ is bounded by
\[
\frac{1}{\min\{m_\eta,1\}}
\int_{\delta_\eta}^{\delta}a(r)^2w(r)\dd r<\infty
\]
uniformly in $\varepsilon\in(0,1)$. We denote this contribution by $O(1)$.
On $(0,\delta_\eta)$, write the integrand of $J_\varepsilon$ as
$N_\varepsilon(r)/D_\varepsilon(r)$, where
\[
N_\varepsilon(r)=(1-\varepsilon)^2a(r)^2w(r),\qquad
D_\varepsilon(r)=(1-\varepsilon)\psi(r)+\varepsilon.
\]
Set
\[
A_\varepsilon=(1-\varepsilon)^2\beta^2\kappa^2g,\qquad
c_\pm=(1-\varepsilon)(1\pm\eta)\kappa.
\]
The preceding bounds on $\psi$, $a$, and $w$ give
\[
(1-\eta)^3A_\varepsilon r^{2\beta-2}
\le N_\varepsilon(r)\le
(1+\eta)^3A_\varepsilon r^{2\beta-2},
\qquad
c_-r^\beta+\varepsilon\le D_\varepsilon(r)
\le c_+r^\beta+\varepsilon.
\]
Therefore, writing
$J^{(\eta)}_\varepsilon=\int_0^{\delta_\eta}
N_\varepsilon(r)/D_\varepsilon(r)\dd r$, we obtain
\begin{equation}\label{eq:squeeze}
(1-\eta)^3A_\varepsilon\!\int_0^{\delta_\eta}\!
\frac{r^{2\beta-2}}{c_+r^\beta+\varepsilon}\dd r
\;\le\;J^{(\eta)}_\varepsilon\;\le\;
(1+\eta)^3A_\varepsilon\!\int_0^{\delta_\eta}\!
\frac{r^{2\beta-2}}{c_-r^\beta+\varepsilon}\dd r.
\end{equation}
Note that
$A_\varepsilon\to\beta^2\kappa^2g$ and
$c_\pm\to(1\pm\eta)\kappa$ as $\varepsilon\downarrow0$.

For $\beta=1$, we use
\[
\int_0^{\delta_\eta}\frac{\dd r}{cr+\varepsilon}
=\frac1c\log\left(1+\frac{c\delta_\eta}{\varepsilon}\right).
\]
Apply this identity to both sides of \eqref{eq:squeeze}. Since
$J_\varepsilon=J^{(\eta)}_\varepsilon+O(1)$, division by
$\log(1/\varepsilon)$ gives
\[
\frac{(1-\eta)^3}{1+\eta}\,\kappa g
\le\liminf_{\varepsilon\downarrow0}
\frac{J_\varepsilon}{\log(1/\varepsilon)}
\le\limsup_{\varepsilon\downarrow0}
\frac{J_\varepsilon}{\log(1/\varepsilon)}
\le\frac{(1+\eta)^3}{1-\eta}\,\kappa g.
\]
Now let $\eta\downarrow0$ to obtain (ii).

For $\tfrac12<\beta<1$, set
$p=(2\beta-1)/\beta$ and $q=(1-\beta)/\beta$. For any $c>0$, the substitution
$r=(\varepsilon/c)^{1/\beta}u$ gives
\[
\int_0^{\delta_\eta}\frac{r^{2\beta-2}}{cr^\beta+\varepsilon}\dd r
=c^{-p}\varepsilon^{-q}
\int_0^{R_\varepsilon}\frac{u^{2\beta-2}}{1+u^\beta}\dd u,
\qquad R_\varepsilon=\delta_\eta(c/\varepsilon)^{1/\beta}\to\infty .
\]
The last integral converges to $I(\beta)$ in \eqref{eq:Ibeta}. Its value there
follows from $v=u^\beta$ and the identity
$\int_0^\infty v^{s-1}/(1+v)\dd v=\pi/\sin(\pi s)$, applied with
$s=2-1/\beta\in(0,1)$.

Apply the change of variables above to both sides of \eqref{eq:squeeze}. Since
$J_\varepsilon=J^{(\eta)}_\varepsilon+O(1)$ and $q>0$, multiplication by
$\varepsilon^q$ gives
\[
\frac{(1-\eta)^3}{(1+\eta)^p}\,B
\le\liminf_{\varepsilon\downarrow0}\varepsilon^qJ_\varepsilon
\le\limsup_{\varepsilon\downarrow0}\varepsilon^qJ_\varepsilon
\le\frac{(1+\eta)^3}{(1-\eta)^p}\,B,
\qquad
B=\beta^2g\kappa^{1/\beta}I(\beta).
\]
Now let $\eta\downarrow0$ to obtain (iii).

For the uniformity claim, the uniform asymptotic assumptions allow
$\delta_\eta$ to be chosen independently of $\alpha$. The contribution of
$[\delta_\eta,\delta)$ is uniformly bounded by the assumed away-from-zero
bounds. In (ii), the logarithmic limit is uniform when the coefficient $c$
ranges in a compact subset of $(0,\infty)$. In (iii), the convergence to
$I(\beta)$ is also uniform in such $c$, because the lower limit of
$R_\varepsilon$ over that compact set tends to infinity. Finally, the
domination used in (i) is uniform near zero and on
$[\delta_\eta,\delta)$. Since $(\kappa_\alpha,g_\alpha)$ ranges in a compact
subset of $(0,\infty)^2$, these observations give the three stated uniform
conclusions.
\end{proof}

\subsection{Proof of Theorem~\ref{thm:main}}

\begin{proof}
Step 1. 
Since $\rho$ is continuous with compact support $\overline D$ and $G_T$ is
smooth and positive, $\phi=\rho/G_T$ is continuous with $\{\phi>0\}=D$, and on
every compact subset of $D$, it is Lipschitz and bounded away from $0$ by
Assumption~\ref{ass:boundary}(b). Moreover, $\nabla\phi=0$ a.e.\ on
$\R^d\setminus\overline D$. From $\nabla\phi=\nabla\rho/G_T-\rho\nabla
G_T/G_T^2$ and $|\nabla G_T|/G_T=|x|/T$, which is bounded on $\overline D$, the
second term is $O(\rho)=O(r^\beta)=o(r^{\beta-1})$. Hence, setting
\[
\kappa(y)=\frac{c(y)}{G_T(y)},\qquad y\in\partial D,
\]
Assumption~\ref{ass:boundary}(c) gives, uniformly in $y\in\partial D$ as
$r\downarrow0$,
\begin{equation}\label{eq:phi-boundary}
\phi\bigl(y+r\nu(y)\bigr)=\kappa(y)r^\beta(1+o(1)),
\qquad
\bigl|\nabla\phi\bigl(y+r\nu(y)\bigr)\bigr|=\beta\kappa(y)r^{\beta-1}(1+o(1)),
\end{equation}
and $G_T(y+r\nu(y))=G_T(y)(1+O(r))$. Since $c$ and $G_T$ are continuous and
positive on the compact set $\partial D$, the pair $(\kappa(y),G_T(y))$ ranges
in a compact subset of $(0,\infty)^2$.

Step 2.
Since $\partial D$ is of class $C^2$ and compact, there is $\delta>0$ such
that $\Psi(y,r):=y+r\nu(y)$ is a $C^1$ diffeomorphism from
$\partial D\times(0,\delta)$ onto the collar
$D_\delta=\{x\in D:\dist(x,\partial D)<\delta\}$, with
$\dist(\Psi(y,r),\partial D)=r$. Let
$S_y=-D\nu(y)|_{T_y\partial D}$ be the shape operator, and denote its
eigenvalues by $\varkappa_1(y),\dots,\varkappa_{d-1}(y)$. Set
\[
J(y,r)=\prod_{i=1}^{d-1}\bigl(1-r\varkappa_i(y)\bigr)=1+O(r)
\quad\text{uniformly in }y.
\]

For $d\ge2$, let $X:U\subset\R^{d-1}\to\partial D$ be a local $C^2$
parametrization and let
\[
G_{ij}(u)=\left\langle\frac{\partial X}{\partial u_i}(u),
\frac{\partial X}{\partial u_j}(u)\right\rangle .
\]
The induced surface measure is
$\dd\sigma(X(u))=\sqrt{\det G(u)}\dd u$. For
$\widetilde\Psi(u,r)=X(u)+r\nu(X(u))$, a direct calculation gives
\[
\bigl|\det D\widetilde\Psi(u,r)\bigr|
=\det\bigl(I-rS_{X(u)}\bigr)\sqrt{\det G(u)}
=J(X(u),r)\sqrt{\det G(u)}.
\]
The usual change-of-variables theorem in each boundary chart, followed by a
partition of unity, therefore gives, for every $F\in L^1(D_\delta)$,
\begin{equation}\label{eq:collar}
\int_{D_\delta}F\dd x
=\int_{\partial D}\int_0^\delta
F\bigl(\Psi(y,r)\bigr)J(y,r)\dd r\dd\sigma(y).
\end{equation}
For $d=1$, \eqref{eq:collar} follows directly from the one-dimensional
change-of-variables formula, with $J\equiv1$ and $\sigma$ the counting measure
on the finite set $\partial D$. The formula also holds for every nonnegative
measurable $F$, by applying the integrable case to $F\wedge m$ and using
monotone convergence.

Take $F=|\nabla\phi|$ in \eqref{eq:collar}. By
\eqref{eq:phi-boundary},
$|\nabla\phi(\Psi(y,r))|\lesssim r^{\beta-1}$, which is integrable because
$\beta>0$. On the compact set $D\setminus D_\delta$, the gradient is bounded
a.e.\ by local Lipschitz regularity. Thus $\phi\in W^{1,1}(D)$ and its continuous
boundary value is zero. The zero-extension theorem for Sobolev functions on
Lipschitz domains then gives
$\phi\in W^{1,1}_{\mathrm{loc}}(\R^d)$ with weak gradient zero a.e.\ outside
$\overline D$. Moreover, $\mu$ is compactly supported and hence has all
exponential moments. Theorem~\ref{thm:reg} therefore applies whenever
$\Ecal_\varepsilon<\infty$.

Step 3.
We split the energy into a boundary part and an interior part. Set
\[
q_\varepsilon(x)
:=\frac{(1-\varepsilon)^2|\nabla\phi(x)|^2}
{(1-\varepsilon)\phi(x)+\varepsilon},
\qquad
R_\varepsilon:=\int_{D\setminus D_\delta}q_\varepsilon\dd\gamma_T.
\]
Since $\nabla\phi=0$ a.e.\ on $\R^d\setminus\overline D$, we have
$\Ecal_\varepsilon=\int_Dq_\varepsilon\dd\gamma_T$.
The interior set $D\setminus D_\delta$ is compact. Hence, by Step 1,
\[
\phi_{\min}:=\min_{x\in D\setminus D_\delta}\phi(x)>0,
\qquad
L:=\mathop{\mathrm{ess\,sup}}_{x\in D\setminus D_\delta}|\nabla\phi(x)|<\infty.
\]
For $0<\varepsilon<1$,
\[
q_\varepsilon
\le \frac{(1-\varepsilon)|\nabla\phi|^2}{\phi}
\le \frac{L^2}{\phi_{\min}}
\quad\text{on }D\setminus D_\delta,
\]
while $q_1=0$. Therefore
$0\le R_\varepsilon\le L^2/\phi_{\min}=O(1)$ uniformly in $\varepsilon$.

We next consider the boundary part. Since $\dd\gamma_T=G_T\dd x$, applying
\eqref{eq:collar} to $F=q_\varepsilon G_T$ gives
\[
\Ecal_\varepsilon
=R_\varepsilon+\int_{\partial D}J_\varepsilon(y)\dd\sigma(y),
\qquad
J_\varepsilon(y)=\int_0^\delta
\frac{(1-\varepsilon)^2\bigl|\nabla\phi(\Psi(y,r))\bigr|^2\,
G_T(\Psi(y,r))\,J(y,r)}{(1-\varepsilon)\phi(\Psi(y,r))+\varepsilon}\dd r. 
\]
For each $y\in\partial D$, this is the quantity in
Lemma~\ref{lem:boundary} with
\[
\psi_y(r)=\phi(\Psi(y,r)),\qquad
a_y(r)=|\nabla\phi(\Psi(y,r))|,\qquad
w_y(r)=G_T(\Psi(y,r))J(y,r),
\]
and parameters $\kappa=\kappa(y)$ and $g=G_T(y)$. Notice that $a_y$ is the
full gradient of $\phi$, rather than the radial derivative $|\psi_y'|$,
because the Fisher energy contains $|\nabla\phi|^2$.

It remains to verify that the lemma applies uniformly in $y$. The three
asymptotic relations near $r=0$ follow uniformly from
\eqref{eq:phi-boundary} and $J(y,r)=1+O(r)$. We may reduce $\delta$ so that the
collar coordinates are defined on a slightly larger interval. Then, for every
$r_0\in(0,\delta)$, the points $\Psi(y,r)$ with
$y\in\partial D$ and $r_0\le r\le\delta$ form a compact subset of $D$.
The local Lipschitz continuity and positivity of $\phi$ on this set give the
two away-from-zero bounds in Lemma~\ref{lem:boundary}, uniformly in $y$.
Hence the conclusions of that lemma hold uniformly on $\partial D$. Finally,
$\sigma(\partial D)<\infty$, so the resulting limits may be integrated over
$y\in\partial D$.

Step 4.
Part (iv) is immediate from Lemma~\ref{lem:boundary}(iv). For (i),
Lemma~\ref{lem:boundary}(i) and the bound on $R_\varepsilon$ give
$\Ecal_\varepsilon\le\Ecal_0<\infty$, the inequality also following directly
from $(1-\varepsilon)^2|\nabla\phi|^2/((1-\varepsilon)\phi+\varepsilon)
\le(1-\varepsilon)|\nabla\phi|^2/\phi$; inserting into \eqref{eq:reg-bound} with
$\varepsilon_n=n^{-1/2}$ gives the stated bound. For (ii),
Lemma~\ref{lem:boundary}(ii) gives
$J_\varepsilon(y)/\log(1/\varepsilon)\to\kappa(y)G_T(y)=c(y)$ uniformly; thus
the factor $G_T$ cancels. Since
$R_\varepsilon/\log(1/\varepsilon)\to0$, the stated limit follows.
Consequently, $\Ecal_\varepsilon\le C_0(1+\log(1/\varepsilon))$ for all
$\varepsilon\in(0,1]$: this follows from the limit for small $\varepsilon$ and,
for $\varepsilon$ bounded away from $0$, from
$\Ecal_\varepsilon\le\varepsilon^{-1}\int|\nabla\phi|^2\dd\gamma_T<\infty$.
Substituting into \eqref{eq:reg-bound} with $\varepsilon_n=n^{-1/2}$ yields
$\delta\le\frac12\sqrt{TC_0(1+\frac12\log n)/n}+n^{-1/2}=O(\sqrt{\log n/n})$.
For (iii), Lemma~\ref{lem:boundary}(iii) gives
$\varepsilon^{(1-\beta)/\beta}J_\varepsilon(y)\to
\beta^2G_T(y)\kappa(y)^{1/\beta}I(\beta)
=\beta^2I(\beta)\,G_T(y)^{1-1/\beta}c(y)^{1/\beta}$ uniformly, and
$\beta^2I(\beta)=\pi\beta/|\sin(\pi/\beta)|$ by \eqref{eq:Ibeta}; integrating
over $\partial D$ and using $R_\varepsilon=O(1)$ gives $\Lambda_\beta$. Hence
$\Ecal_\varepsilon\le C\varepsilon^{-q}$ with $q=(1-\beta)/\beta$ for all
$\varepsilon\in(0,1]$. Enlarging $C$ to cover $\varepsilon$ bounded away from
$0$ and applying Remark~\ref{rem:q} gives the rate
$n^{-1/(q+2)}=n^{-\beta/(1+\beta)}$.
\end{proof}

\subsection{Triangular kernel densities}

We now apply Theorem~\ref{thm:main}(ii) to finite mixtures of triangular kernel densities.

\begin{corollary}\label{cor:triangular}
Let $d=1$, $m\in\mathbb N$, $h_0>0$, $x_1,\dots,x_m\in\R$ and
\[
\rho=\sum_{j=1}^m\rho_j,\qquad
\rho_j(x)=\frac1{mh_0}\Bigl(1-\Bigl|\frac{x-x_j}{h_0}\Bigr|\Bigr)_+ .
\]
Then Assumption~\ref{ass:boundary} holds with $\beta=1$, and if $q_k^-$ (resp.\
$q_k^+$) denotes the number of components whose supports have left endpoint
$a_k$ (resp.\ right endpoint $b_k$), then, with
$\varepsilon_n:=n^{-1/2}$,
\[
\lim_{\varepsilon\downarrow0}\frac{\Ecal_\varepsilon}{\log(1/\varepsilon)}
=\frac1{mh_0^2}\sum_{k=1}^K\bigl(q_k^-+q_k^+\bigr),
\qquad
\delta\bigl(\Law(Y_n^{\varepsilon_n}),\mu\bigr)
=O\Bigl(\sqrt{\tfrac{\log n}{n}}\Bigr).
\]
\end{corollary}

\begin{proof}
Each $\rho_j$ is supported on $[x_j-h_0,x_j+h_0]$, is piecewise linear with
slopes $\pm1/(mh_0^2)$, and satisfies $\rho_j(x_j\mp h_0\pm r)=r/(mh_0^2)$ for
$0<r<h_0$. Hence $\{\rho>0\}$ is a finite union of bounded open intervals and
$\rho$ is continuous and piecewise linear, so
Assumption~\ref{ass:boundary}(a) and (b) hold. At a left endpoint $a_k$, exactly the
$q_k^-$ components with left endpoint $a_k$ contribute, giving
$\rho(a_k+r)=q_k^-r/(mh_0^2)$ for small $r>0$, i.e.\ $c_k^-=q_k^-/(mh_0^2)$, and
likewise $c_k^+=q_k^+/(mh_0^2)$. Apply Theorem~\ref{thm:main}(ii) together with
$\int_{\partial D}c\dd\sigma=\sum_k(c_k^-+c_k^+)$.
\end{proof}

For this example, \cite{EndoNakano2024} established qualitative two-stage
convergence. Corollary~\ref{cor:triangular} provides the admissible coupling
$\varepsilon_n=n^{-1/2}$ and the corresponding rate.

\begin{remark}
Compared with Proposition~\ref{prop:generic}, the gain comes from the size of
the integrand near the boundary. With $r=\dist(x,\partial D)$ one has
$\phi\asymp r$ and
$|\nabla\phi|\asymp1$ when $\beta=1$, so
$|\nabla\phi|^2/(\phi+\varepsilon)\asymp(r+\varepsilon)^{-1}$, whose integral
grows like $\log(1/\varepsilon)$ rather than $\varepsilon^{-1}$. A slower growth
of $\Ecal_\varepsilon$ permits a smaller $\varepsilon$, hence a smaller
regularization bias.
\end{remark}

\section{Concluding remarks and future directions}\label{sec:remarks}

Theorem~\ref{thm:main} gives the sharp growth rates and leading constants of
the regularized Fisher energy. Combined with the entropy bound and Pinsker's
inequality, these asymptotics yield the exponents $\frac12$, $\frac12$ up to a
logarithmic factor, and $\frac\beta{1+\beta}$. A natural next step is to seek
matching lower bounds for $\delta(\Law(Y_n^{\varepsilon_n}),\mu)$ and to compare
the mixture regularization \eqref{eq:reg-intro} with other schemes.

For $\beta\le\frac12$, Theorem~\ref{thm:main}(iv) shows that
$\Ecal_\varepsilon=\infty$. Since \eqref{eq:reg-intro} makes $\phi$ positive but
leaves $\nabla\phi$ unchanged, this regime calls for regularizations that smooth
both quantities, for example by convolving $\rho$ at a scale coupled to $n$ and
balancing the regularization and discretization errors.

The $C^2$ assumption on $\partial D$ is used only for the collar formula
\eqref{eq:collar}. This suggests extending the result to piecewise smooth
domains by applying the argument on smooth faces and controlling neighborhoods
of corners and edges. Since these singular parts have zero surface measure, the
leading constants are expected to remain unchanged.

\subsection*{Acknowledgments}

This work was supported by JSPS KAKENHI Grant Number JP24K06861.

\end{document}